\documentclass{conm-p-l}
\usepackage{amssymb,amsmath}
\usepackage[mathscr]{eucal}

\newcommand{\al}{\alpha}
\newcommand{\eu}{\stackrel{?}{=}}
\DeclareMathOperator*{\GL}{\mathrm{GL}}
\DeclareMathOperator*{\Int}{\int}
\newcommand{\la}{\lambda}
\DeclareMathOperator{\Li}{\mathrm{Li}}
\newcommand{\ou}{\overline{1}}
\newcommand{\Q}{\mathbf{Q}}
\newcommand\shuff{
  \setlength{\unitlength}{.4pt}
  \begin{picture}(40,20)
    \put(10,2){\line(1,0){20}} \put(10,2){\line(0,1){10}}
    \put(20,2){\line(0,1){10}} \put(30,2){\line(0,1){10}}
  \end{picture}}
\newcommand\qsh{\setlength{\unitlength}{.4pt}
  \begin{picture}(40,20)
    \put(10,2){\line(1,0){20}} \put(10,2){\line(0,1){10}}
    \put(20,2){\line(0,1){10}} \put(30,2){\line(0,1){10}}
  \end{picture}_{\!q}\;}
\newcommand{\Sh}{{\mathrm{Sh}}}
\newcommand{\Shuff}{{\mathrm{Shuff}}}
\newcommand{\e}{\eta}
\newcommand{\w}{\omega}
\newcommand{\Shq}{{\mathrm{Shq}}}
\DeclareMathOperator{\sinc}{\mathrm{sinc}}
\newcommand{\z}{\zeta}
\newcommand{\R}{\mathbf{R}}
\newcommand{\C}{\mathbf{C}}
\newcommand{\Z}{\mathbf{Z}}

\newtheorem{theorem}{Theorem}[section]
\newtheorem{Cor}[theorem]{Corollary}
\newtheorem{prop}[theorem]{Proposition}
\newtheorem{lemma}[theorem]{Lemma}

\theoremstyle{definition}
\newtheorem{definition}[theorem]{Definition}
\newtheorem{example}[theorem]{Example}

\theoremstyle{remark}

\numberwithin{equation}{section}

\begin{document}

\title{Multiple Polylogarithms: A Brief Survey}

\author{Douglas Bowman}
\address{University of Illinois at Urbana-Champaign,
Department of Mathematics,
273 Altgeld Hall, 1409 W. Green St., Urbana, IL 61801 U.S.A.}
\email{bowman@math.uiuc.edu}
\thanks{The first author was supported in part by NSF Grant DMS-9705782.}

\author{David M.~Bradley}
\address{Department of Mathematics and Statistics,
University of Maine, 5752 Neville Hall, Orono, ME 04469-5752
U.S.A.} \email{dbradley@e-math.ams.org}
\thanks{The second author was
partially supported by the University of Maine summer
faculty research fund.}

\copyrightinfo{2000}{American Mathematical Society}
\subjclass{Primary 33E20; Secondary 11G55, 11M99, 40B05}
\date{January 28, 2001.}


\keywords{Euler sums, multiple zeta values,
          polylogarithms, multiple harmonic series,
          quantum field theory, knot theory, Riemann zeta function}

\begin{abstract}
We survey various results and conjectures concerning multiple
polylogarithms and the multiple zeta function. Among the results,
we announce our resolution of several conjectures on multiple zeta
values.  We also provide a new integral representation for the
general multiple polylogarithm, and develop a $q$-analogue of the
shuffle product.

\end{abstract}

\maketitle

\specialsection{Introduction}

In recent years, nested harmonic sums have attracted increasing
attention in both the mathematics and physics communities.   The
sums occur within the context of knot theory and quantum field
theory, yet their rich structure offers much to fascinate
theoreticians in such diverse areas as algebra, number theory, and
combinatorics.  Multiple polylogarithms generalize the
aforementioned nested sums, as well as the Riemann zeta function
and the classical polylogarithm, while still retaining many
interesting properties. Their study has led to many beautiful yet
unproven conjectures, including evaluations at arbitrary depth
discovered with the use of recently developed integer
relations-finding algorithms and high precision numerical
computations in the thousands of digits.

Multiple polylogarithms~\cite{BBBLa,Gonch1,Gonch3} are multiply
nested sums of the form
\begin{equation}
   {\Li}_{s_1,\dots,s_k}(z_1,\dots,z_k)
   :=\sum_{n_1>\cdots>n_k>0}\;\prod_{j=1}^k n_j^{-s_j}
   {z_j}^{n_j},
\label{Li-nest}
\end{equation}
where $s_1,\dots,s_k$ and $z_1,\dots,z_k$ are complex numbers
suitably restricted so that the sum~(\ref{Li-nest}) converges.
Instances of multiple polylogarithms have occurred in several
disparate fields, such as combinatorics (analysis of
quad-trees~\cite{FLLS,LL} and of lattice reduction
algorithms~\cite{DFV}), knot
theory~\cite{BGK,BK2,BK1,LeM,LeM2,Tak}, perturbative quantum field
theory~\cite{BarBroad,DJB0,DJB00,DJB1} and mirror
symmetry~\cite{Hoff4}. There is also quite sophisticated work
relating polylogarithms and their generalizations to arithmetic
and algebraic geometry, and to algebraic
$K$-theory~\cite{Beil,Brow2,Brow,Gonch1,Gonch2,Gonch3,Woj3,Woj2,Woj}.

Figuring prominently are the nested sums~(\ref{Li-nest}) in which
each $z_j=1$.  These latter are now commonly referred to as
multiple zeta values~\cite{BBBLc,BBBLa,BowBrad2,YOhno1,Zag}
and
are denoted by
\begin{equation}
\label{MZVdef}
   \zeta(s_1,\dots,s_k) := \sum_{n_1>\cdots>n_k>0}\;
   \prod_{j=1}^k n_j^{-s_j}.
\end{equation}
The study of such sums goes back to Euler~\cite{LE}, who showed
that
\begin{equation}
   2\z(m,1) = m\z(m+1)-\sum_{j=1}^{m-2}\z(m-j)\z(j+1),
   \qquad 2\le m\in\Z.
\label{Euler}
\end{equation}
It can be shown that the sum~(\ref{MZVdef}) is absolutely
convergent in the region
\[
   \{(s_1,\dots,s_k)\in\C^k : \sum_{j=1}^r \Re(s_j)> r
   \;\mbox{for}\; 1\le r\le k\}.
\]
(The condition given in Proposition 1 of~\cite{Zhao} is
insufficient to guarantee absolute convergence.)

Define the {\it depth} of the multiple
polylogarithm~(\ref{Li-nest}) to be the number $k$ of nested
summations.  A good deal of work on multiple polylogarithms, and
more specifically multiple zeta values, has been motivated by the
problem of determining which sums can be expressed (say
polynomially with rational coefficients) in terms of other sums of
lesser depth. Settling this question in complete generality is
currently beyond the reach of number theory.  For example, proving
the irrationality of expressions such as
$\zeta(5,3)/\zeta(5)\zeta(3)$ appears to be impossible with
current techniques. Nevertheless, considerable progress has been
made with regard to proving specific classes of reductions, even
at arbitrary depth.  The first nontrivial success at arbitrary
depth was the settling~\cite{BBBLc,BBBLa} of Zagier's
conjecture~\cite{Zag}
\begin{equation}
   \zeta(\underbrace{3,1,3,1,\ldots,3,1}_{2n})
   =4^{-n}\zeta(\underbrace{4,4,\ldots,4}_{n})
   = \frac{2\pi^{4n}}{(4n+2)!},
   \qquad 0\le n\in\Z,
\label{Z31}
\end{equation}
in which the $2n$ and $n$ beneath the underbraces in~(\ref{Z31})
denote the depth of the respective multiple zeta values.
Subsequent work (see eg.~\cite{BowBrad1,BowBrad2,BowBrad3})
has
largely focused on developing a suitable framework for dealing
with ultimately periodic argument strings in general, and
additionally sums of multiple zeta values whose set of argument
strings is fixed by the action of certain subgroups of the group
of permutations.

It is instructive to trace the development of the subject and see
for oneself how ad hoc techniques and considerations have in many
cases evolved into more systematic methods of ongoing interest. In
this connection, one might begin by citing the partial fractions
technique of Euler~\cite{LE} and Nielsen~\cite{Niels},
subsequently employed by many others
eg.~\cite{BBG,BG,Hoff1,HoffCoen,Mark,Mord,SubSit}, and which Ohno
recently parlayed in his exceedingly clever proof of the cyclic
sum formula~\cite{HoffOhno,YOhno2}.  Techniques based on
elementary integration formul{\ae} and identities for special
functions tailored to specific examples eg.~\cite{BB,BBG,Ded,Mord}
have evolved~\cite{BBB,BBBLa} into quite general, sophisticated,
and powerful analytic methods~\cite{BowBrad1,BowBrad3}.  The
na\"{\i}ve approach of deriving elementary series transformation
identities and solving the resulting systems of linear
equations~\cite{BG,RaoSub}, used to prove reducibility results of
depth three or less, has been largely superceded (eg.\ by methods
based on contour integration~\cite{Flaj}) and supplanted by
considerations of the shuffle and stuffle~\cite{BBBLa}
multiplications, and relatedly the harmonic algebra and the
algebra of quasi-symmetric functions~\cite{Hoff2,Hoff3,HoffOhno}.

Computational issues---both numeric and symbolic---have also come
into play.  Relations satisfied by multiple polylogarithms, and
multiple zeta values in particular, can be exploited by symbolic
computer algebra systems to prove reductions of small
weight~\cite{MinhPet1}. (Here the {\em weight} of the multiple
zeta value~(\ref{MZVdef}) is simply the sum of the arguments
$s_1+\cdots+s_k$.)  Interest in high-precision, rapid computation
of multiple zeta values~\cite{Cran,Cran-Buh} (see
also~\cite[\S7.2]{BBBLa}) has been stimulated by the ability to
numerically hunt for or rule out identities (to a high degree of
probability) with the aid of recently developed integer
relations-finding algorithms~\cite{Forcade,Hastad,Lenstra}.

In addition to the as yet unsolved problem of classifying all
possible relationships between multiple zeta values at positive
integer arguments, one can also consider~(\ref{MZVdef}) as a
function of the complex numbers $s_1,\dots,s_k$ and consider
questions regarding analytic continuation, trivial zeros, and
values at the non-positive integers. The analytic continuation
of~(\ref{MZVdef}) in the case $k=2$ was established by
Atkinson~\cite{Atk} via the Poisson summation formula, and later
by Apostol and Vu~\cite{Apost}, who used the Euler-Maclaurin
summation formula. Subsequently, Arakawa and Kaneko~\cite{AK}
proved that if $s_2,\dots,s_k$ are fixed positive integers,
then~(\ref{MZVdef}) can be meromorphically continued as a function
of $s_1$ to the whole complex $s_1$-plane. The analytic
continuation of~(\ref{MZVdef}) as a function defined on $\C^k$ was
established by Akiyama, Egami and Tanigawa~\cite{Aki1} using the
Euler-Maclaurin summation formula. An independent approach due to
Zhao~\cite{Zhao} uses properties of Gelfand and Shilov's
generalized functions~\cite{GS}.  Zhao also attempts a discussion
of trivial zeros for $k\le 3$.
To our knowledge,
no-one has yet determined the trivial zeros of~(\ref{MZVdef}) for
general $k$.

The issue of values of~(\ref{MZVdef}) at the non-positive integers
is subtle, since for $k>1$ the result will in general depend on
the order in which the respective limits are taken.  Thus, for
example, if $n$ is a non-negative integer, $s(k,j)$ and $S(k,j)$
denote the Stirling numbers of the first and second kind,
respectively, and $B_j$ denotes the $j$th Bernoulli number,
then~\cite{Aki2}
\[
   \lim_{s_1\to-n}\lim_{s_2\to0}\cdots\lim_{s_k\to0}
   \z(s_1,\dots,s_k) = \frac{(-1)^{n+1}}{n+1}
   \sum_{j=1}^{n+1} \frac{(-1)^{k+j} j! S(n+1,j)}{k+j},
\]
whereas
\[
   \lim_{s_k\to0}\cdots\lim_{s_2\to0}\lim_{s_1\to -n}
   \z(s_1,\dots,s_k) = (-1)^k\delta_{n,0}
   -\frac{1}{(k-1)!}\sum_{j=1}^k
   \frac{s(k,j)B_{n+j}}{n+j},
\]
where $\delta_{n,0}=0$ if $n>0$ and $\delta_{0,0}=1$.  In
particular (cf.\ also~\cite{Zhao})
\[
   \lim_{s_1\to 0}\lim_{s_2\to 0} \z(s_1,s_2) = \frac{1}{3},
   \quad \mbox{but}\quad
   \lim_{s_2\to 0}\lim_{s_1\to 0} \z(s_1,s_2) = \frac{5}{12}.
\]
We are unaware of any systematic treatment in the case of
\emph{arbitrary} non-positive integer arguments.

\subsection{Notation}
Let $\sigma_1,\dots,\sigma_k\in\{-1,1\}$.   We will have occasion
to discuss the particular sums of the form
\begin{equation}
    \Li_{s_1,\dots,s_k}(\sigma_1x,\sigma_2,\dots,\sigma_k)
    = \sum_{n_1>\cdots>n_k>0} x^{n_1}
   \prod_{j=1}^k n_j^{-s_j} \sigma_j^{n_j},
\label{MZVxDef}
\end{equation}
in which $0\le x\le 1$ is real and $s_1,\dots,s_k$ are positive
integers with $x=s_1=\sigma_1=1$ excluded for convergence.
Accepted practice dictates that~(\ref{MZVxDef}) may be abbreviated
by $\z_x(s_1,\dots,s_k)$ with a bar placed over $s_j$ if and only
if $\sigma_j=-1$.  When $x=1$, these are called Euler sums.  Thus
a multiple zeta value is an Euler sum with no alternations. We
adopt the convention that $\z_x() := 1$ when no arguments are
present ($k=0$).  We also drop the subscript $x$ when $x=1$ since
$\z_1(s_1,\dots,s_k)$ agrees with~(\ref{MZVdef}) when each $s_j$
is bar-free.  For example,
\[
   \z(\overline{2},1) = \sum_{n=1}^\infty \frac{(-1)^n}{n^2}
   \sum_{m=1}^{n-1}\frac{1}{m}.
\]
It will be convenient to abbreviate strings of repeated arguments
by using exponentiation to denote concatenation powers. Then the
first two members of~(\ref{Z31}) may be written
$\z(\{3,1\}^n)=4^{-n}\z(\{4\}^n)$.

Finally, as customary the Gaussian hypergeometric function and the
logarithmic derivative of the Euler gamma function  are denoted by
\[
   F(a,b;c;x) = \sum_{n=0}^\infty x^n \prod_{j=0}^{n-1}
   \frac{(a+j)(b+j)}{(1+j)(c+j)}
   \quad{\text{and}}\quad
   \psi=\frac{\Gamma'}{\Gamma},
\]
respectively.  We also abbreviate the set of the first $k$
positive integers $\{1,2,\dots,k\}$ by $N_k$.

\specialsection{Stuffles}

For the sake of brevity and simplicity, we shall restrict the
discussion in this section to multiple zeta values.  For a
discussion of the more general polylogarithmic case,
see~\cite{BBBLa}.

As Hoffman~\cite{HoffOhno} observed, one can view multiple zeta
values as values of a homomorphism on a commutative $\Q$-algebra
in two ways; the $\Q$-algebra multiplications have been referred
to elsewhere~\cite{BBBLa} as ``shuffle'' and ``stuffle.''  It is
conjectured that all relations between multiple zeta values are a
consequence of the collision of the two multiplications, provided
one admits the divergent sums~(\ref{MZVdef}) with $s_1=1$
(suitably renormalized) into the model. However, there seems
little hope of proving this conjecture in the near future, and at
present a wide variety of analytic, algebraic, and combinatorial
techniques are used to prove identities for multiple zeta values.

Stuffle relations, or more simply stuffles---see
\S\ref{sect:DrinInt} and \S\ref{sect:cyclic} below for a
discussion of shuffles---arise when one multiplies two nested
series of the form~(\ref{MZVdef}) and expands the product
distributively.  Thus if $u$ and $v$ are (ordered) lists of
positive integers, then
\[
   \z(u)\z(v) = \sum_{w\in u*v} \z(w),
\]
where $u*v$ is the multiset defined by the recursion
\begin{equation}
   su * tv = s(u* tv) \cup t( su * v) \cup (s+t)(s*t),
   \qquad 1\le s,t\in\Z.
\label{StuffleDef}
\end{equation}
In~(\ref{StuffleDef}) it is to be understood that if $M$ is a
multiset of lists and $a$ is an integer, then $aM$ denotes the
multiset of lists obtained by placing $a$ at the front of each
list in $M$. For example $(s,t)*u =
\{(s,t,u),(s,t+u),(s,u,t),(s+u,t),(u,s,t)\}$ and correspondingly
$\z(s,t)\z(u)=\z(s,t,u)+\z(s,t+u)+\z(s,u,t)+\z(s+u,t)+\z(u,s,t).$

Let $f(|u|,|v|)$ denote the number of lists in $u*v$.  The
recursive decomposition~(\ref{StuffleDef}) shows that the
generating function
\[
   F(x,y) := \sum_{m=0}^\infty \sum_{n=0}^\infty f(m,n) x^m y^n
\]
satisfies the functional equation
$F(x,y)=1+xF(x,y)+yF(x,y)+xyF(x,y)$.  It follows that
$F(x,y)=(1-x-y-xy)^{-1}$ and hence that
\begin{equation}
   f(m,n) = \sum_{k=0}^m \binom{m}{k}\binom{n+k}{m}
          = \sum_{k=0}^{\min(m,n)}\binom{n}{k}\binom{m}{k} 2^k.
\label{NumStuff}
\end{equation}
One can also give a direct, combinatorial proof
of~(\ref{NumStuff}) by considering how the indices interlace in
the product of two nested series of the form~(\ref{MZVdef}).

There are interesting connections between stuffles, polyominoes,
and codes which we briefly indicate.  To begin, note that a
stuffle counted by $f(m,n)$ can be viewed as a pair $(\phi,\psi)$
of order-preserving injections
\[
   \phi:N_m\to N_r,\qquad \psi:N_n\to N_r
\]
where $r$ is chosen so that $\max(m,n)\le r\le m+n$ and the union
of the images of $\phi$ and $\psi$ is all of $N_r$.  One can
associate to such a pair a sequence of integers $b_1,\dots,b_m$ by
defining $a_1=\phi(1)-1$ and $a_j=\phi(j)-\phi(j-1)-1$ for $2\le
j\le m$ and then letting
\[
   b_j= \left\{\begin{array}{ll}
    -a_j &\mbox{if  $\phi(j)$ is in the image of $\psi$},\\
    a_j &\mbox{otherwise}\end{array}\right.
\]
for each $j\in N_m$. Since $\phi$ is order-preserving, $a_j\ge 0$
for each $j\in N_m$ and $\sum_{j=1}^m |b_j| = \phi(m)-m\le n$.
Conversely, given a sequence of integers $b_1,\dots,b_m$
satisfying $\sum_{j=1}^m |b_j|\le n$, the pair $(\phi,\psi)$ is
uniquely determined.  Let $p=|\{j:b_j<0\}|$.   We have $r=m+n-p$,
$\phi(1)=|b_1|+1$ and $\phi(j)=\phi(j-1)+|b_j|+1$ for $2\le j\le
m$.  Put $\psi(j)=\phi(j)$ if $b_j<0$.  The remaining values of
$\psi$ are determined by the requirement that it be an
order-preserving injective map of $N_n$ to $N_r$.  Thus, there is
a one-to-one correspondence between the stuffles counted by
$f(n,m)$ and the sets of integer lattice points whose
cardinalities satisfy
\begin{equation}
   \bigg|\bigg\{(b_1,\dots,b_m)\in \Z^m : \sum_{j=1}^m |b_j|\le n\bigg\}
   \bigg|
   = \bigg|\bigg\{(b_1,\dots,b_n)\in \Z^n : \sum_{j=1}^n |b_j|\le
   m\bigg\}\bigg|,
\label{LatticeSets}
\end{equation}
the identity~(\ref{LatticeSets}) holding in view of the obvious
symmetry $f(m,n)=f(n,m)$.

Define an $n$-dimensional polyomino formed by adding $m$ coats to
a single-celled polyomino, where a coat consists of just enough
cells to cover each previously exposed $(n-1)$-dimensional face.
There is clearly a bijection between such polyominos and the
second set of lattice points~(\ref{LatticeSets}).  The
relationship is explored in greater detail in~\cite{GW}.

\specialsection{Integral Representations}

\subsection{The Drinfeld Integral}\label{sect:DrinInt}
There is also a representation for multiple zeta values in terms
of an ``iterated'' (Drinfeld) integral due to
Kontsevich~\cite{Zag}. For real $0\le x\le 1$ and positive
integers $s_1,\dots,s_k$ with $x=s_1=1$ excluded for convergence,
we have
\begin{equation}
   \z_x(s_1,\dots,s_k) =\int \prod_{j=1}^k
   \bigg(\prod_{r=1}^{s_j-1} \frac{dt_r^{(j)}}{t_r^{(j)}}\bigg)
   \frac{dt_{s_j}^{(j)}}{1-t_{s_j}^{(j)}},
\label{iterint}
\end{equation}
where the integral is over the simplex
\[
x>t_1^{(1)}>\cdots>t_{s_1}^{(1)}>\cdots>t_1^{(k)}>\cdots>t_{s_k}^{(k
)}>0,
\]
and is abbreviated by
\begin{equation}
   \int_0^x \prod_{j=1}^k a^{s_j-1}b,
   \qquad a=dt/t,\quad b=dt/(1-t).
\label{shortiterint}
\end{equation}
Making the simultaneous change of variable $t\mapsto 1-t$ at each
level of integration and then reversing the order of integration
makes transparent the duality identity for multiple zeta values:
\begin{equation}
   \z(s_1+2,\{1\}^{r_1},\dots,s_k+2,\{1\}^{r_k})
   = \z(r_k+2,\{1\}^{s_k},\dots,r_1+2,\{1\}^{s_1}),
\label{duality}
\end{equation}
first conjectured in~\cite{Hoff1} and proved in~\cite{Zag}.

A related integral representation enabled Ohno~\cite{YOhno1} to
prove the following beautiful generalization of~(\ref{duality}).
Let
\[
   S(p_1,\dots,p_n;m) := \sum_{c_1+\cdots+c_n=m}
   \z(p_1+c_1,\dots,p_n+c_n),
\]
where the sum is over all non-negative integers $c_1,\dots,c_n$
which sum to $m$.  As in~(\ref{duality}) define the dual argument
lists
\[
   p :=(s_1+2,\{1\}^{r_1},\dots,s_k+2,\{1\}^{r_k})
\]
and
\[
   p' :=(r_k+2,\{1\}^{s_k},\dots,r_1+2,\{1\}^{s_1}).
\]
Then~\cite{YOhno1} $S(p;m) = S(p';m).$  When $m=0$, Ohno's
result
reduces to~(\ref{duality}).  Another interesting specialization is
obtained by taking $p=(k+1)$ and $m=n-k-1$; one then deduces
Granville's theorem~\cite{Granv}, originally conjectured
independently by Courtney Moen~\cite{Hoff1} and Michael
Schmidt~\cite{Mark}:
\[
   \sum_{s_1+\cdots+s_k=n} \z(s_1,\dots,s_k) = \z(n),
\]
where the sum is over all positive integers $s_1,\dots,s_k$ which
sum to $n$ and $s_1>1$.

The iterated integral representation is also responsible for a
second multiplication rule satisfied by multiple zeta values.
Suppose that $x,y\in\R$ and $f_j:[y,x]\to\R$ are integrable
functions for $j=1,2,\dots,n$.  It is customary to make the
abbreviation
\begin{equation}
   \int_y^x \prod_{j=1}^n \alpha_j :=
   \Int_{x>t_1>t_2>\cdots>t_n>y} \;\prod_{j=1}^n f_j(t_j)\,dt_j,
   \qquad \alpha_j := f_j(t_j)\,dt_j,
\label{IterIntNotn}
\end{equation}
with the convention that~(\ref{IterIntNotn}) is equal to 1 if
$n=0$ regardless of the values of $x$ and $y$.
There is an alternative definition of iterated integrals which
explains their name.  For $j=1,2,\dots,n$ again define
the $1$-forms $\al_j$ by $\al_j := f_j(t_j)\,dt_j$.
Then put
\begin{eqnarray}\label{itdef2}
   \int_y^x \al_1\al_2\cdots\al_n
   &:=& \left\{ \begin{array}{ll}
   \int_y^x f_1(t_1) \int_y^{t_1}\al_2\cdots\al_n \,dt_1 &\mbox{if
$n>0$}\\
   1 &\mbox{if $n=0$.}\\
   \end{array}\right.
\end{eqnarray}
Expanding out this second definition, it is easy to see that
it coincides with the definiton as an integral over a simplex.
Both definitions occur frequently in the literature.

Clearly the product
of two iterated integrals of the form~(\ref{IterIntNotn}) consists
of a sum of iterated integrals involving all possible interlacings
of the variables. Therefore, if we denote the set of all
$(m+n)!/m!\,n!$ permutations $\sigma$ of the indices $N_{m+n}$
satisfying $\sigma^{-1}(j)<\sigma^{-1}(k)$ for all $1\le j<k\le m$
and $m+1\le j<k\le m+n$ by $\Shuff(m,n)$, then we have the
self-evident formula
\begin{equation}\label{shufprdmotiv}
   \bigg(\int_y^x \prod_{j=1}^m \alpha_j \bigg)\bigg(\int_y^x\,
   \prod_{j=m+1}^{m+n}\alpha_j\bigg)
   = \sum_{\sigma\in\Shuff(m,n)} \int_y^x\; \prod_{j=1}^{m+n}
   \alpha_{\sigma(j)},
\end{equation}
and so define the shuffle product $\shuff$ by
\begin{equation}
\label{shuff-rule}
   \bigg(\prod_{j=1}^m \alpha_j\bigg)\shuff
   \bigg(\prod_{j=m+1}^{m+n}\alpha_j\bigg)
    := \sum_{\sigma\in\Shuff(m,n)} \prod_{j=1}^{m+n}
    \alpha_{\sigma(j)}.
\end{equation}
Thus, the sum is over all non-commutative products (counting
multiplicity) of length $m+n$ in which the relative orders of the
factors in the products $\alpha_1\alpha_2\cdots \alpha_m$ and
$\alpha_{m+1}\alpha_{n+2}\cdots \alpha_{m+n}$ are preserved. The
term ``shuffle'' is used because such permutations arise in riffle
shuffling a deck of $m+n$ cards cut into one pile of $m$ cards and
a second pile of $n$ cards.

The study of shuffles and iterated integrals was pioneered by
Chen~\cite{Chen54,Chen57} and subsequently formalized by
Ree~\cite{Ree}.  As with the case of stuffles, one can view an
element of $\Shuff(m,n)$ as a pair of order-preserving injections
$(\phi,\psi)$ where now $\phi:N_m\to N_{m+n}$ and $\psi:N_n\to
N_{m+n}$ have disjoint images.  One can then define a vector
$(a_1,\dots,a_m)$ of non-negative integers by $a_1=\phi(1)-1$ and
$a_j=\phi(j)-\phi(j-1)-1$ for $2\le j\le m$.  Since $\phi$ is
order-preserving, $a_j\ge 0$ for each $j\in N_m$ and $\sum_{j=1}^m
a_j=\phi(m)-m\le n.$ Conversely, if we have such a vector of
non-negative integers, then $\phi(1)=a_1+1$ and
$\phi(j)=\phi(j-1)+a_j+1$ for $2\le j\le m$ defines an
order-preserving injection $\phi:N_m\to N_{m+n}$, and hence a
shuffle.  Thus, there is a one-to-one correspondence between
$\Shuff(m,n)$  and the sets of non-negative integer lattice points
whose cardinalities satisfy
\[
   \bigg|\bigg\{(a_1,\dots,a_m) \in \Z_{\ge 0}^m: \sum_{j=1}^m a_j\le
    n\bigg\}\bigg|
    = \bigg| \bigg\{(a_1,\dots,a_n) \in \Z_{\ge 0}^n: \sum_{j=1}^n a_j\le
    m\bigg\}\bigg|,
\]
the latter identity holding in light of the fact that
$\Shuff(m,n)$ is clearly symmetric in $m$ and $n$.  A deeper study
of the algebra and combinatorics of shuffles leads to an
alternative proof of~(\ref{Z31}) and generalizations thereof; see
\S\ref{sect:cyclic}.

\subsection{A New Integral Representation}
In light of the usefulness of the various integral
representations, it may be of interest to give here a new integral
representation for~(\ref{Li-nest}).  The new representation
appears to embody properties of both the Drinfeld and partition
integrals of~\cite{BBBLa}, and therefore may be useful in proving
certain results for multiple polylogarithms which have thus far
withstood attacks based on traditional methods.  The derivation
employs MacMahon's Omega operator, which discards terms with
non-positive exponents from formal Laurent series in
$\la_1,\dots,\la_k$. Thus, in view of~(\ref{Li-nest}), if $0\le
x_1,\dots,x_k\le 1$, we may write
\begin{eqnarray*}
   &&\Li_{s_1,\dots,s_k}(x_1,\dots,x_k)\\
   &=& \Omega \prod_{j=1}^k \sum_{n_j>0} n_j^{-s_j}
   \left(x_j\la_j\la_{j-1}^{-1}\right)^{n_j},\qquad \la_0:=1\\
   &=& \Omega \prod_{j=1}^k
   \Li_{s_j}\left(x_j\la_j\la_{j-1}^{-1}\right)\\
   &=& \Omega \prod_{j=1}^k
   \Int_{1>u_1^{(j)}>\cdots>u_{s_j}^{(j)}>0} \bigg(\prod_{r=1}^{s_j-1}
      \frac{du_r^{(j)}}{u_r^{(j)}}\bigg)\frac{x_j\la_j\la_{j-1}^{-1}
      \,du_{s_j}^{(j)}}{1-x_j\la_j\la_{j-1}^{-1}u_{s_j}^{(j)}}\\
   &=& \Omega \prod_{j=1}^k
   \Int_{1>u_1^{(j)}>\cdots>u_{s_j}^{(j)}>0} \bigg(\prod_{r=1}^{s_j-1}
      \frac{du_r^{(j)}}{u_r^{(j)}}\bigg)
      \sum_{m_j=1}^\infty
      \left(x_j\la_j\la_{j-1}^{-1}\right)^{m_j}
      \big(u_{s_j}^{(j)}\big)^{m_j-1}\,du_{s_j}^{(j)}\\
   &=& \int_{\Delta(\vec s)}\left\{\prod_{j=1}^k\bigg(\prod_{r=1}^{s_j-
1}
   \frac{du_r^{(j)}}{u_r^{(j)}}\bigg)\right\}
  \sum_{m_1>\cdots>m_k>0}\;\prod_{j=1}^k
   \left(x_ju_{s_j}^{(j)}\right)^{m_j}\frac{du_{s_j}^{(j)}}{u_{s_j}^{(j)}},
\end{eqnarray*}
where $\Delta(\vec s)$ denotes the set of all integration
variables satisfying
\[
   1>u_1^{(j)}>u_2^{(j)}>\cdots>u_{s_j}^{(j)}>0
\]
for $j=1,2,\dots,k$.  Since $0\le y_j<1$ for each $j=1,2,\dots,
k$ implies
\begin{eqnarray*}
   \sum_{m_1>\cdots>m_k>0}\; \prod_{j=1}^k y_j^{m_j}
   &=& \sum_{n_1=1}^\infty\cdots \sum_{n_k=1}^\infty
       y_1^{n_1+\cdots+n_k}y_2^{n_2+\cdots+n_k}\cdots y_k^{n_k}\\
   &=& \frac{y_1}{1-y_1}\cdot\frac{y_1y_2}{1-y_1y_2}\cdots
       \frac{y_1y_2\cdots y_k}{1-y_1y_2\cdots y_k},
\end{eqnarray*}
it follows that
\begin{equation}
   \Li_{s_1,\dots,s_k}(x_1,\dots,x_k)
   = \int_{\Delta(\vec s)} \prod_{j=1}^k \left\{
       \tau\bigg(\prod_{m=1}^j x_m u_{s_m}^{(m)}\bigg)
      \prod_{r=1}^{s_j}\frac{du_r^{(j)}}{u_r^{(j)}}\right\},
\end{equation}
where $\tau(x):=x/(1-x).$

\specialsection{Generating Functions}

In many cases, generating functions provide the best means of
stating reductions involving one or more parameters.  A specific
example of this which also illustrates how knowledge of the
subject has progressed is given first.  We then outline a
systematic approach for tackling multiple zeta values with
periodic argument lists, followed by additional examples to
illustrate the richness of the theory.

\subsection{Two-Parameter Symmetry}
In connection with Euler's result~(\ref{Euler}),
Markett~\cite{Mark} derived
\begin{eqnarray}
   \z(s,1,1) &=&
   \frac{1}{6}s(s+1)\z(s+2)-\frac{1}{2}(s-1)\z(2)\z(s)
   -\frac{s}{4}\sum_{n=0}^{s-4}\z(s-n-1)\z(n+3)\nonumber\\
   &+&\frac{1}{6}\sum_{n=0}^{s-4}\z(s-n-2)\sum_{m=0}^n\z(n-
m+2)\z(m+2),
   \qquad 3\le s\in\Z,
\label{Mark}
\end{eqnarray}
via elementary but intricate series manipulations and partial
fraction identities.  An equivalent formula is proved
in~\cite{BBG} using elementary facts about the dilogarithm, the
polygamma function and the higher derivatives of the Euler beta
function.
%
For larger values of $n$, the representation of
$\z(s,\{1\}^n)$ in terms of values of the Riemann zeta function
becomes increasingly complicated.  Nevertheless, there is an
elegant generating function formulation which we restate here.
\begin{theorem}[\cite{BBB}]\label{thm:Drin}
The bivariate formal power series identity
\begin{multline}
   \sum_{m=0}^\infty \sum_{n=0}^\infty
   x^{m+1}y^{n+1}\z(m+2,\{1\}^n)\\
   =1-\exp\bigg\{\sum_{k=2}^\infty
   \frac1{k}\left(x^k+y^k-(x+y)^k\right)\z(k)\bigg\}
\label{Drin}
\end{multline}
holds.
\end{theorem}

\begin{Cor} Let $n$ and $s$ be non-negative integers with $s\ge
2$.  Then $\z(s,\{1\}^n)$ lies in the polynomial ring
$\Q[\z(2),\z(3),\dots,\z(s+n)]$.
\end{Cor}

By comparing coefficients of $x^{s-1}y^{n+1}$ on both sides
of~(\ref{Drin}), one sees that in fact, $\z(s,\{1\}^n)$ is a
rational linear combination of products of Riemann zeta values
such that the sum of the arguments in each product is equal to
$s+n$.
%
Moreover, Euler's result~(\ref{Euler}) is an immediate consequence
of comparing coefficients of $x^{s-1}y^2$.  Similarly, Markett's
formula~(\ref{Mark}) can be obtained most easily by comparing
coefficients of $x^{s-1}y^3$. Finally, as the right hand side
of~(\ref{Drin}) is evidently symmetric in $x$ and $y$, the left
hand side must also be.  Thus Theorem~\ref{thm:Drin} implies the
special case $\z(m+2,\{1\}^n)=\z(n+2,\{1\}^m)$ of the duality
formula~(\ref{duality}).  It would be interesting to find a
generating function formulation of duality at full strength.

\subsection{Periodic Argument Lists}
Results such as~(\ref{Z31}) and~(\ref{Drin}) suggest that one
might profit from a more systematic study of multiple zeta values
whose argument lists form an ultimately periodic sequence.  This
is indeed the case; such a study forms the basis of some of our
current work in progress~\cite{BowBrad3}.

\subsubsection{Period One}

The case of all identical arguments is quite well understood.
Nevertheless, there are a few items of interest worth recording
here, in particular a connection to the problem of determining the
number of unordered factorizations of an integer.

For $\Re(s)>1$, equation~(\ref{MZVdef}) implies
\begin{equation}
   \sum_{k=0}^\infty t^{ks}\z(\{s\}^k)
   = \prod_{j=1}^\infty \bigg(1+\frac{t^s}{j^s}\bigg).
\label{Period1GF}
\end{equation}
If in~(\ref{Period1GF}) we take $s$ to be an even integer, say
$s=2n$ where $n$ is a positive integer, then we may
rewrite~(\ref{Period1GF}) in the form
\begin{equation}
   \sum_{k=0}^\infty (-1)^k t^{2kn} \z(\{2n\}^k)
   =\prod_{j=0}^{n-1}\sinc(\pi t \rho^j),
\label{sincs}
\end{equation}
where $\rho = e^{\pi i/n}$ and $\sinc x = \sin x/x$ for $x\ne 0$;
$\sinc 0:=1$.  The identity~(\ref{sincs}) is one of many possible
generalizations of Euler's formula for $\z(2n)$, and moreover
shows that $\z(\{2n\}^k)$ is a rational multiple of $\pi^{2kn}$.

Differentiating both sides of~(\ref{Period1GF}) and equating
coefficients yields the recurrence
\begin{equation}
   k\z(\{s\}^k) = \sum_{j=1}^k (-1)^{j+1}
   \z(js)\z(\{s\}^{k-j}),
   \qquad 0\le k\in\Z,\quad \Re(s)>1,
\label{CrandallREC}
\end{equation}
which is really just a special case of Newton's formula
\[
   ke_k = \sum_{j=1}^k (-1)^{j+1}p_je_{k-j},\qquad 0\le k\in\Z,
\]
relating the elementary symmetric functions and power sum
symmetric functions
\[
   e_k = \sum_{j_1>\cdots>j_k>0} \;\prod_{m=1}^k x_{j_m},
   \qquad
   p_k := \sum_{j>0} x_j^k.
\]
Substituting $1/j^s$ for each indeterminate $x_j$ yields $e_k =
\z(\{s\}^k)$ and $p_k=\z(ks)$.

From~(\ref{CrandallREC}) it follows that if $k$ is a positive
integer and $\Re(s)>1$, then $\z(\{s\}^k)$ lies in the polynomial
ring $\Q[\z(s),\z(2s),\dots,\z(ks)]$.  In fact, there is an
explicit formula for $\z(\{s\}^k)$ in terms of a sum over
partitions of $k$.

\begin{definition}\label{PrtnDef}
Let $r$ be a non-negative integer and let
$\al=(\al_1,\al_2,\dots)$ be a non-negative integer partition of
$r$.  Let $m_j=\#\{i:\al_i=j\}$ be the number of parts of size
$j$, and put $c_{\al} = \prod_{j\ge 1} m_j!(-j)^{m_j}.$
Furthermore, abbreviate $r=\sum_{j\ge 1}\al_j$ by $|\al|$ and
$\prod_{j\ge 1}p_{\al_j}$ by $p_{\al}$.
\end{definition}

In view of the generic relationship~\cite{IJMac}
\[
   \sum_{k=0}^\infty e_k t^k = \exp\bigg\{-\sum_{r=1}^\infty
   \frac{(-t)^r p_r}{r}\bigg\}
   = \sum_{\al} (-t)^{|\al|} c_{\al}^{-1} p_{\al},
\]
for $\Re(s)>1$ we therefore have
\[
   \sum_{k=0}^\infty t^k \z(\{s\}^k) =
   \exp\bigg\{-\sum_{r=1}^\infty \frac{(-1)^r \z(rs)t^r}{r}\bigg\}
   =\sum_{\al} (-t)^{|\al|} c_{\al}^{-1} \prod_{\al_j>0}
   \z(\al_j s),
\]
i.e.\
\begin{equation}
   \z(\{s\}^k) = (-1)^k \sum_{|\al|=k}
   c_{\al}^{-1}\prod_{\al_j>0}\z(\al_j s).
\label{Period1Explicit}
\end{equation}

We note the following connection with factorisatio
numerorum~\cite{Hensley}. (See
also~\cite{Canfield,Oppenheim,Warlimont}.)
Let $\al$ be as in Definition~\ref{PrtnDef}.  Define the
unrestricted divisor function associated with the partition $\al$
by
\[
   d_{\al}(m) =\sum_{\prod_{j\ge 1}d_j^{\al_j}=m}1.
\]
For example $d_{1,1}$ is the ordinary divisor function, and
$d_2(m)=1$ if $m$ is a perfect square and zero otherwise.
\begin{prop} Let $\tau_k(m)$ denote the number of unordered
factorizations of $m$ into $k$ distinct factors.  Then
\[
   \tau_k(m)=(-1)^k \sum_{|\al|=k} c_{\al}^{-1}d_{\al}(m).
\]
\end{prop}
\begin{proof}
Observe that for $\Re(s)>1$,
\[
   \z(\{s\}^k) = \sum_{n_1>\cdots>n_k>0}\;\prod_{j=1}^k n_j^{-s}
   = \sum_{m=1}^\infty \tau_k(m)\, m^{-s}.
\]
Now compare coefficients of $m^{-s}$ in~(\ref{Period1Explicit}).
\end{proof}
\begin{example} Since $\z(\{s\}^2) = \tfrac12 \z^2(s)-\tfrac12\z(2s)$,
we get $\tau_2(m)=\tfrac12 d_{1,1}(m)-\tfrac12 d_2(m)$.  In
particular $\tau_2(12)= 3$.
\end{example}

\subsubsection{Period Two and Beyond}

In contrast with the situation in which all arguments are
identical, much remains to be explored in the case of argument
strings of period two and higher.  In~\cite{BBBLa}
and~\cite{BowBrad1} differential equations were found to be a
useful technique for analyzing the generating functions for period
2. We summarize here some results from~\cite{BBBLa}
and~\cite{BowBrad1} to indicate the richness and complexity of the
resulting formul{\ae} arising from the solution of the associated
fourth order differential equation.
\begin{definition}\label{hyperdefs} For $0\le x\le 1$ and
$z\in\C$, let
\begin{eqnarray*}
   Y_1(x,z) &:=& F(z,-z;1;x),\\ 
   Y_2(x,z) &:=& (1-x)F(1+z,1-z;2;1-x),\\ 
   G(z) &:=&
   \tfrac14\left\{\psi(1+iz)+\psi(1-iz)-\psi(1+z)-\psi(1-z)\right\}.
\end{eqnarray*}
\end{definition}

\begin{theorem}[\cite{BBBLa}]\label{thm:zx31gf}
Let $Y_1$ be as in Definition~\ref{hyperdefs}.  Then for $0\le
x\le 1$ and $|z|<1$,
\begin{equation}
   \sum_{n=0}^\infty (-1)^n z^{4n}
   4^n\z_x(\{3,1\}^n) = Y_1(x,z)Y_1(x,iz).
\label{ZFact}
\end{equation}
\end{theorem}

\begin{theorem}[\cite{BowBrad1}]\label{prop:zx313gf}
Let $Y_1$, $Y_2$ and $G$ be as
in Definition~\ref{hyperdefs}.  Then for $0\le x\le 1$ and
$|z|<1$,
\begin{multline}
   \sum_{n=0}^\infty (-1)^n z^{4n+2} 4^n \z_x(3,\{1,3\}^n)
   =G(z)Y_1(x,z)Y_1(x,iz)\\-\frac{Y_1(x,iz)Y_2(x,z)}{4Y_1(1,z)}
   +\frac{Y_1(x,z)Y_2(x,iz)}{4Y_1(1,iz)}.
\label{z313GF}
\end{multline}
\end{theorem}
Note that~(\ref{ZFact}) proves~(\ref{Z31}).
Similarly~(\ref{z313GF}) proves
\begin{eqnarray*}
   \z(3,\{1,3\}^n) &=& 4^{-n}\sum_{k=0}^n\z(4k+3)\z(\{4\}^{n-k})\\
   &=& \sum_{k=0}^n
   \frac{2\pi^{4k}}{(4k+2)!}\left(-\frac{1}{4}\right)^{n-k}\z(4n-4k+3),
\end{eqnarray*}
which escaped the extensive numerical and symbolic searches
carried out in the preparation of~\cite{BBB,BBBLc,BBBLa}.
Differentiation of~(\ref{z313GF}) followed by a delicate analysis
of the asymptotic behaviour of the requisite hypergeometric
functions at their singular points proves~\cite{BowBrad1} the
reduction
\begin{multline*}
   \z(2,\{1,3\}^n)\\
    = 4^{-n}\sum_{k=0}^n(-1)^k\z(\{4\}^{n-k})
   \bigg\{(4k+1)\z(4k+2)-4\sum_{j=1}^k\z(4j-1)\z(4k-4j+3)\bigg\}
\end{multline*}
conjectured in~\cite{BBB,BBBLa}.

The proof of Theorem~\ref{prop:zx313gf} hinges on showing that
both sides of~(\ref{ZFact}) and~(\ref{z313GF}) are annihilated by
the same fourth order differential operator.  In~\cite{BBBLa},
computer algebra was used to establish this for~(\ref{ZFact}).  At
the time, a conceptual proof was unavailable. Subsequently the
present authors (see~\cite{BowBrad1}) found a conceptual proof of
the following more general result, which is perhaps best
understood in the context of work going back to Orr~\cite{Orr} and
Clausen~\cite{Claus} on differential equations satisfied by a
product of hypergeometric series.  The result is shown
in~\cite{BowBrad2} to be closely related to the combinatorial
``shuffle'' approach outlined in \S\ref{sect:cyclic}, and may be
stated as follows.
\begin{lemma}\label{lem:zfac} Let $K$ be a differential field of
characteristic not equal to $2$ and let $D$ be a derivation on
$K$. For each $k\in K,$ define a derivation $D_k:=kD$.  Let $t$ be
a constant, and suppose that for some $f,g,u,v\in K$ the
differential equations $(D_fD_g+t)u=0$ and $(D_fD_g-t)v=0$ hold.
Then $uv$ is annhilated by the differential operator
$(D_f^2D_g^2+4t^2).$
\end{lemma}
In particular, taking $f(x)=1-x$, $g(x)=x$ and $t=z^2$, given that
$Y_1$ and $Y_2$ satisfy $(D_fD_g+z^2)y=0$,
Lemma~\ref{lem:zfac}
shows that each of the three linearly independent functions
$Y_1(x,z)Y_1(x,iz)$, $Y_1(x,iz)Y_2(x,z)$, and $Y_1(x,z)Y_2(x,iz)$
are annihilated by the operator $D_f^2D_g^2+4z^4$.  That $L(x,z)$
and $S(x,z)$ are annihilated by the same operator follows easily
from the integral representation~(\ref{iterint}), whence
Theorem~\ref{prop:zx313gf} is proved.

Since $D_f^2D_g^2+4z^4$ is a fourth order differential operator,
one might legitimately ask in what context the fourth linearly
independent solution $Y_2(x,z)Y_2(x,iz)$ arises.  It turns out
that due to the double logarithmic singularity arising from the
product of the underlying hypergeometric functions at $x=1$, it is
easier to ascribe a meaning to this solution in the case of
alternating sums~(\ref{MZVxDef}).  Recalling the generating
function
\begin{equation}
   A(z):=\sum_{n=0}^{\infty}z^n\z(\{\ou\}^n) =
   \prod_{j=1}^\infty \bigg(1+\frac{(-1)^jz}{j}\bigg)
   =\frac{\Gamma(1/2)}{\Gamma(1+z/2)\Gamma(1/2-z/2)}
\label{Adef}
\end{equation}
from~\cite{BBB}, we have
\begin{theorem}[\cite{BowBrad1}]\label{prop:Mgf}
Let $0\le x\le 1$, and $|t|< \infty $. Put $z=(1+i)t/2$,
$s=(1+x)/2$, and let $U(s,z) = Y_1(s,z)-zY_2(s,z)$, where $Y_1$
and $Y_2$ are as in Definition~\ref{hyperdefs}.  Then,
\begin{equation}
    \sum_{n=0}^\infty
    \bigg[t^{2n}\z_x(\{\ou,1\}^n)+t^{2n+1}\z_x(\ou,\{1,\ou\}^n)\bigg]
     = \frac{U(s,-z)U(s,iz)}{A(-z)A(iz)}.
\label{Mrep}
\end{equation}
\end{theorem}
Theorem~\ref{prop:Mgf} is a bivariate generalization of the
conjecture~\cite[equation~(14)]{BBB} in the case $x=1$, and may be
viewed as an analytic extension of the purely combinatorial
identity~(\ref{MFact}) below.

In recent work~\cite{BowBrad3}, the authors have greatly extended
the differential equation approach. The authors have obtained
results on more general generating functions which include not
only multiple zeta values, but polylogarithmic and
hyperlogarithmic~\cite{LD} values in general. In fact, from the
point of view of iterated integrals, arbitrary forms may occur in
the iterated integrals studied. The differential equations are
still present. The authors have classified various bases for the
solutions of the differential equations, given matrices for change
of basis, and found the explicit representations of the monodromy
matrices of the associated differential equations. These results
actually stand out with greater distinction in a more general
setting. Taking arbitrary $1$-forms on a manifold $M$, an explicit
homomorphism is obtained from $\pi_1(M,x_0)$ into $\GL_n(\C)$.
This gives rise to a transport between the manifold $M$ and its
principle bundle constructed from the representation into
$\GL_n(\C)$.  Finally these results can be cast yet more generally
in the setting of differentiable spaces. Our homomorphism is
similar to the celebrated homomorphism of
K.~T.~Chen~\cite{Chen54,Chen57,Chen58,Chen71} in that it is built
out of a generating function of iterated integrals. The essential
difference is that Chen's homomorphism maps into a formal Lie
group, while our homomorphism maps into $\GL_n(\C)$. Will our
homomorphism give different information than Chen's? We are
currently investigating the geometric implications of our work in
this area.

\specialsection{Shuffles and Cyclic Insertion}\label{sect:cyclic}

As in~\cite{MinhPet1} (cf.\ also \cite{BBBLc,Ree}) let $A$ be a
finite set and let $A^*$ denote the free monoid generated by $A$.
We regard $A$ as an alphabet, and the elements of $A^*$ as words
formed by concatenating any finite number of letters (repetitions
permitted) from the alphabet $A$. By linearly extending the
concatenation product to the set $\Q\langle A\rangle$ of rational
linear combinations of elements of $A^*$, we obtain a
non-commutative polynomial ring with the elements
of $A$ being indeterminates and with
multiplicative identity $1$ denoting the empty
word.

The shuffle product~(\ref{shuff-rule}) is alternatively defined
first on words by the recursion
\begin{equation}
\begin{cases}
   \forall w\in A^*, \quad & 1\shuff w = w\shuff 1 = w,\\
   \forall a,b\in A, \quad\forall u,v\in A^*, \quad & au\shuff bv
   =a(u\shuff bv)+b(au\shuff v),
\end{cases}
\label{LeftShuffDef}
\end{equation}
and then extended linearly to $\Q\langle A\rangle$.  One checks
that the shuffle product so defined is associative and
commutative, and thus $\Q\langle A\rangle$ equipped with the
shuffle product becomes a commutative $\Q$-algebra, denoted
$\Sh_{\Q}[A]$.  Radford~\cite{Rad} has shown that $\Sh_{\Q}[A]$ is
isomorphic to the polynomial algebra $\Q[L]$ obtained by adjoining
the transcendence basis $L$ of Lyndon words to the field $\Q$ of
rational numbers.

The recursive definition~(\ref{LeftShuffDef}) has its analytical
motivation in the formula for integration by parts---equivalently,
the product rule for differentiation.  Thus, if we put
$a=f(t)\,dt$, $b=g(t)\,dt$ and
\[
   F(x) := \int_y^x (au\shuff bv) =
   \bigg(\int_y^x f(t)\int_y^t u\,dt\bigg)\bigg(\int_y^x g(t)\int_y^t
   v\,dt\bigg),
\]
then writing $F(x)=\int_y^x F'(s)\,ds$ and applying the product
rule for differentiation yields
\begin{eqnarray*}
   F(x) &=& \int_y^x \bigg(f(s)\int_y^s u\bigg)\bigg(\int_y^s
   g(t)\int_y^t v\,dt\bigg)\,ds\\
   &&\qquad\qquad\qquad + \int_y^x g(s)\bigg(\int_y^s
   f(t)\int_y^t u\,dt\bigg)\int_y^s v\,ds\\
   &=& \int_y^x \left[ a(u\shuff bv) + b(au\shuff v)\right].
\end{eqnarray*}
Alternatively, by viewing $F$ as a function of $y$, we see that
the recursion~(\ref{LeftShuffDef}) could equally well have been
stated as
\begin{equation}
\begin{cases}
   \forall w\in A^*, \quad & 1\shuff w = w\shuff 1 = w,\\
   \forall a,b\in A, \quad\forall u,v\in A^*, \quad & ua\shuff vb
   =(u\shuff vb)a+(ua\shuff v)b.
\end{cases}
\label{RightShuffDef}
\end{equation}
Of course, both definitions are equivalent to~(\ref{shuff-rule}).

The combinatorial proof~\cite{BBBLc} of Zagier's
conjecture~(\ref{Z31}) hinged on expressing the sum of the words
comprising the shuffle product of $(ab)^n$ with $(ab)^m$ as a
linear combination of basis subsums. In~\cite{BowBrad2} a more
comprehensive study of the shuffle algebra $\Sh_{\Q}[a,b]$ is
undertaken, and as a consequence correspondingly deeper results
for multiple zeta values are obtained.  To highlight the most
interesting of these results, we first recall the following
\begin{definition}[\cite{BBBLc}]\label{def:T}
For integers $m\ge n\ge 0$ let $S_{m,n}$ denote the set of words
occurring in the shuffle product $(ab)^{n}\shuff (ab)^{m-n}$ in
which the subword $a^2$ appears exactly $n$ times, and let
$T_{m,n}$ be the sum of the $m!/(2n)!(m-2n)!$ distinct words in
$S_{m,n}.$  For all other integer pairs $(m,n)$ it is convenient
to define $T_{m,n}:=0$.
\end{definition}
One then has
\begin{theorem}[\cite{BowBrad2}]\label{thm:T-Binom}
Let $x$ and $y$ be commuting indeterminates, and let $m$ be a
non-negative integer. In the commutative polynomial ring
$(\Sh_{\Q}[a,b])[x,y]$ we have the shuffle convolution formula
\begin{equation}
\label{T-Binom}
   \sum_{k=0}^m x^k y^{m-k} \left[(ab)^k \shuff (ab)^{m-k}\right]
   = \sum_{n=0}^{\lfloor m/2\rfloor} (4xy)^n (x+y)^{m-2n}\,
   T_{m,n}.
\end{equation}
\end{theorem}
A special case of Theorem~\ref{thm:T-Binom} implies the intriguing
shuffle factorization due to Broadhurst, and which in turn
implies~(\ref{Z31}):
\begin{equation}
   A\bigg(\frac{z}{1-i}\bigg) \shuff\, A\bigg(\frac{z}{1+i}\bigg)
   = M(z) \in (\Sh_{\Q}[a,b])[[z]],\quad i^2=-1,
\label{MFact}
\end{equation}
where
\[
    A(z) := \sum_{n=0}^\infty (z^2ab)^n(1+za)
    \quad {\text{and}}\quad
    M(z):= \sum_{n=0}^\infty (z^4a^2b^2)^n(1+za+z^2a^2+z^3a^2b).
\]
The experts will recognize~(\ref{Adef}) as the analytic version of
$A(z)$ above, in which $a=-dt/(1+t)$ and $b=dt/(1-t)$.  Similarly
for $M(z)$ and the left hand size of~(\ref{Mrep}) when $x=1$.

In addition, Theorem~\ref{thm:T-Binom} plays a key role in a
remarkable combinatorial generalization of~(\ref{Z31}) which we
proceed to describe.  Let $S_{m,n}$ be as in
Definition~\ref{def:T}.  Note that each word in $S_{m,n}$ has a
unique representation
\begin{equation}
   (ab)^{m_0}\prod_{k=1}^n (a^2 b)(ab)^{m_{2k-1}}b(ab)^{m_{2k}},
\label{CompositionCorrespondence}
\end{equation}
in which $m_0,m_1,\dots,m_{2n}$ are non-negative integers with
sum
$m-2n$.  Conversely, every ordered $(2n+1)$-tuple
$(m_0,m_1,\dots,m_{2n})$ of non-negative integers with sum $m-2n$
gives rise to a unique word in $S_{m,n}$
via~(\ref{CompositionCorrespondence}).  Thus, a bijective
correspondence $\varphi$ is established between the set $S_{m,n}$
and the set $C_{2n+1}(m-2n)$ of ordered non-negative integer
compositions of $m-2n$ with $2n+1$ parts.  In view of the
relationship~(\ref{iterint}) expressing multiple zeta values as
iterated integrals, it therefore makes sense to define
\[
   Z(\vec s) := \int_0^1 \varphi(\vec s),\qquad \vec s\in C_{2n+1}(m-
2n),
\]
where as in~(\ref{shortiterint}), we now identify the abstract
letters $a$ and $b$ with the differential 1-forms $dt/t$ and
$dt/(1-t)$, respectively.  Thus, if $\vec
s=(m_0,m_1,\dots,m_{2n})$, then
\begin{eqnarray*}
   Z(\vec s) &=&
   \int_0^1 (ab)^{m_0}\prod_{k=1}^n (a^2b)(ab)^{m_{2k-
1}}b(ab)^{m_{2k}}\\
   &=&\z(\{2\}^{m_0},3,\{2\}^{m_1},1,\{2\}^{m_2},3,\{2\}^{m_3},1,
   \dots,3,\{2\}^{m_{2n-1}},1,\{2\}^{m_{2n}}),
\end{eqnarray*}
in which the argument string consisting of $m_j$ consecutive twos
is inserted after the $j$th element of the string $\{3,1\}^n$ for
each $j=0,1,2,\dots,2n$.  It turns out~\cite{BowBrad2} that
\begin{equation}
    \sum_{\vec s\in C_{2n+1}(m-2n)} Z(\vec s) =
\frac{2\pi^{2m}}{(2m+2)!}
   \binom{m+1}{2n+1},
\label{Zsum}
\end{equation}
for all non-negative integers $m$ and $n$ with $m\ge 2n$.  The
proof uses Theorem~\ref{thm:T-Binom} at essentially full strength
combined with some tricky generatingfunctionology. Observe that
equation~(\ref{Z31}) is the special case of~(\ref{Zsum}) in which
$m=2n$, since $Z(\{0\}^{2n+1})=\z(\{3,1\}^n)$.

A more compelling formulation of~(\ref{Zsum}) can be given as
follows.  Again, let $\vec s=(m_0,m_1,\dots,m_{2n})$ and put
\[
   {\mathscr{C}}(\vec s) :=Z(\vec s)+\sum_{j=1}^{2n}
   Z(m_j,m_{j+1},\dots,m_{2n},m_0,\dots,m_{j-1}).
\]
In other words, sum over all cyclic permutations of the argument
list  $\vec s$.   Then~\cite{BowBrad2}
\begin{equation}
   \sum_{\vec s\in C_{2n+1}(m-2n)}{\mathscr{C}}(\vec s) =
   Z(m)\times|C_{2n+1}(m-2n)| =
   \frac{\pi^{2m}}{(2m+1)!}\binom{m}{2n}
\label{Csum}
\end{equation}
is an equivalent formulation of~(\ref{Zsum}).  Here, we have used
\[
   Z(m) = \z(\{2\}^m) = \frac{\pi^{2m}}{(2m+1)!},
   \qquad 0\le m\in\Z,
\]
which follows from~(\ref{sincs}).  The cyclic insertion
conjecture~\cite{BBBLa} can be restated as the assertion that
${\mathscr{C}}(\vec s)=Z(m)$ for all $\vec s\in C_{2n+1}(m-2n)$
and integers $m\ge 2n\ge 0$. Thus, (\ref{Csum}) reduces the
problem to that of establishing the invariance of
${\mathscr{C}}(\vec s)$ on $C_{2n+1}(m-2n)$.  It is likely that
this remaining step can be accomplished using only the shuffle
property of multiple zeta values in conjunction
with~(\ref{sincs}).

\specialsection{Dimension Conjectures}

Broadhurst~\cite{DJB3} has conjectures concerning the size of
various bases (graded by weight and depth) for expressing multiple
zeta values in terms of either irreducible multiple zeta values,
or irreducible Euler sums, and also for expressing Euler sums in
terms of irreducible Euler sums. The adjunction of additional
differential forms appears to simplify the problem at each stage.
Thus, if $D(n,k)$ denotes the number of multiple zeta values of
weight $n$ and depth $k$ in a minimal $\Q$-basis for reducing all
multiple zeta values to a $\Q$-linear combination of products of
basis multiple zeta values, it is conjectured that
\[
   \prod_{n\ge 3}\prod_{k\ge 1}(1-x^ny^k)^{D(n,k)} \eu
   1-\frac{x^3y}{1-x^2}+\frac{x^{12}y^2(1-y^2)}{(1-x^4)(1-x^6)}.
\]
However, if we allow Euler sums into the basis, letting $M(n,k)$
denote the minimal number of Euler sums of weight $n$ and depth
$k$ needed to reduce all multiple zeta values to basis Euler sums,
then
\[
   \prod_{n\ge 3}\prod_{k\ge 1}(1-x^ny^k)^{E(n,k)}\eu
   1-\frac{x^3y}{(1-x^2)(1-xy)}.
\]
One can also consider the problem of reducing Euler sums in terms
of basis Euler sums.  Let $E(n,k)$ denote the minimal number of
Euler sums of weight $n$ and depth $k$ required to reduce all
Euler sums to basis Euler sums.  It is conjectured that
\[
   \prod_{n\ge 3}\prod_{k\ge 1}(1-x^ny^k)^{E(n,k)} \eu
   1-\frac{x^3y}{1-x^2}.
\]
Adjoining forms associated with sixth roots of unity to the set
of possible differential
forms yields the multiple Clausen values~\cite{BBK}, and here it
is conjectured that the number $P(n,k)$ of irreducible multiple
Clausen values of weight $n$ and depth $k$ is generated by
\[
   \prod_{n>1}\prod_{k>0}(1-x^ny^k)^{P(n,k)} \eu
   1-\frac{x^2y}{1-x}.
\]

\specialsection{$q$-Shuffles}

Here we consider a $q$-analogue of the shuffle algebra discussed
in \S\ref{sect:cyclic}. Let $A$ be a set, not necessarily finite,
and let $\e:A\to A$ be bijective. Now form the free monoid
generated by $A$ and call it $A^*$ as before. Extend the action of
$\e$ to $A^*$ in the obvious way so that $\e$ becomes an
automorphism of $A^*$. Again regard $A$ as an alphabet, and the
elements of $A^*$ as words formed by concatenating any finite
number of letters (repetitions permitted) from the alphabet $A$.
By linearly extending the concatenation product to the set
$\Q\langle A\rangle$ of rational linear combinations of elements
of $A^*$, we obtain a non-commutative polynomial ring with the
elements of $A$ being indeterminates and with multiplicative
identity $1$ denoting the empty word. It is clear that $\e$ now
extends to an automorphism of $\Q\langle A\rangle$.

A $q$-shuffle algebra is defined to be the ordered pair
$(\Q\langle A\rangle,\qsh)$ where $\qsh$ is a commutative and
associative bilinear operator on $\Q\langle A\rangle$ satisfying
the identity
\begin{equation}
\begin{cases}
   \forall w\in A^*, \quad & 1\qsh w = w\qsh 1 = w,\\
   \forall a,b\in A, \quad\forall u,v\in A^*, \quad & au\qsh bv
   =a(u\qsh bv)+b(\e(au)\qsh v).
\end{cases}
\label{q-LeftShuffDef}
\end{equation}
We denote a $q$-shuffle algebra over $A$ by $\Shq_{\Q}[A]$. It
will be observed that a $q$-shuffle algebra is a commutative $\Q$-algebra.

Our definition implies that there may be more than one way of
writing the $q$-shuffle product of two words. For example, letting
$a,b\in A$, it is easy to see that in $\Shq_{\Q}[A]$ one has
$a\qsh b=ab+b\e a=ba+a\e b$. As the length of the words being
multiplied increases the number of different expressions also
grows.

The motivation for our definition of $\Shq_{\Q}[A]$ is not
difficult to see. As the shuffle algebra is motivated by the
property (\ref{shufprdmotiv}) of iterated integrals, one wants a
similar identity to hold for iterated Jackson
$q$-integrals~\cite{Gasp}. Recall the definition of a Jackson
$q$-integral. For $x>0$, let $f:[0,x]\to\R$ be Riemann integrable.
The Jackson $q$-integral of $f$ on $[0,x]$ is defined by
\begin{equation}
   \int_0^x f(t)\,d_qt:=\sum_{n\ge 0}f(xq^n)\,xq^n(1-q).
\label{qintegral}
\end{equation}
Because for any $0<q<1$ the sum on the right hand side
of~(\ref{qintegral}) is a Riemann sum for $\int_0^x f(t)\,dt,$ it
follows that the Jackson $q$-integral tends to the ordinary
Riemann integral in the limit as $q$ approaches 1.

One defines iterated Jackson $q$-integrals in exactly the same way
that ordinary iterated integrals are defined by (\ref{itdef2}). To
this end, for $j=1,2,\dots,n$ let $f_j:[0,x]\to\R$ and $\w_j :=
f_j(t_j)\,d_qt_j$. Then put
\begin{eqnarray}
\label{qiterint}
   \int_0^x \w_1\w_2\cdots\w_n
   &:=& \prod_{j=1}^n\int_0^{t_{j-1}}f_j(t_j)\,d_qt_j,\quad t_0:=x
   \\
   &=& \left\{ \begin{array}{ll}
   \int_0^x f_1(t_1) \int_0^{t_1}\w_2\cdots\w_n \,d_qt_1 &\mbox{if
         $n>0$}\nonumber\\
   1 &\mbox{if $n=0$.}
   \end{array}\right.
\end{eqnarray}
Here the fact that the $1$-forms on the right hand side
of~(\ref{qiterint}) are $q$-difference $1$-forms implies that the
integral on the left hand side of~(\ref{qiterint}) is a
$q$-iterated integral and not an ordinary iterated integral.

Iterating the definition of the $q$-iterated integral, one finds
that
\begin{multline}\label{ugly}
\int_0^x\w_1\cdots\w_k= \sum_{n_1,\dots,n_k\ge 0}
f_1(xq^{n_k})f_2(xq^{n_k+n_{k-1}})\cdots f_k(xq^{n_k+\cdots+n_1})
\\ \times q^{kn_k+(k-1)n_{k-1}+\cdots+2n_2+n_1}(1-q)^kx^k,
\end{multline}
but this is not a very convenient expression. To
simplify~(\ref{ugly}) it helps to make the following change of
indices:
\[
l_i:=\sum_{j=k-i+1}^k n_j, \qquad 1\le i\le k.
\]
Then (\ref{ugly}) reduces to
\begin{multline}
\int_0^x\w_1\cdots\w_k=\\ (1-q)^k\sum_{0\le l_1\le l_2\le\cdots\le
l_k} xq^{l_1}f(xq^{l_1})xq^{l_2}f(xq^{l_2})\cdots
xq^{l_k}f(xq^{l_k}).
\end{multline}
This in turn motivates the definition of $q$-difference $1$-forms
by the equation
\begin{equation}
\w_i=\w_i(t)=f_i(t)\,d_qt:=f_i(t)t(1-q).
\end{equation}
Notice that with this definition, the $q$-integral reduces to a
summation operator on these $1$-forms:
\[
\int_0^x \w = \sum_{j\ge 0} \w(t_0q^j),
\]
which agrees with the original definition by virtue of $t_0=x$.

Now fix $0<q<1$ and define the Rogers~\cite{Rog} $q$-difference operator
$\e$ acting on continuous functions $f:[0,\infty]\to\R$ by the equation
$(\e f)(x):=f(xq)$. Also define the $q$-derivative in the usual
way:
\[
(D_qf)(x):=\frac{f(x)-f(xq)}{x(1-q)}.
\]
The relevant fact about the $q$-derivative is the following
well-known property:
\begin{equation}\label{duh}
D_q\int_0^x f(t)\,d_qt=f(x).
\end{equation}

We are ready to give the motivation for the $q$-shuffle product.
The alphabet $A$ now consists of $q$-difference $1$-forms, and the
automorphism $\e$ is Rogers' $q$-difference operator. The action
of $\e$ is extended to forms by the equation
$\e\w=\w(tq)=f(tq)tq(1-q)$. This of course defines a new form
$\w'= \e\w$ by $\w'=g(t)t(1-q)$, where $g(t)=qf(tq)$. The alphabet
$A$ needs to be infinite to account for all the forms $\e^j\w$ for
$j\in\Z$. The action of $\e$ is now extended to $\Q\langle
A\rangle$ in the obvious way. It clearly forms an automorphism of
this algebra.  We wish to define the $q$-shuffle product so that
for $u,v\in A^*$ the following equation is true:
\begin{equation}\label{motive}
\int_0^x u\qsh v = \bigg(\int_0^x u\bigg)\bigg(\int_0^x v\bigg).
\end{equation}
To accomplish this, one applies the $q$-analogue of the argument
given in \S\ref{sect:cyclic} for deriving the recursive definition
of the shuffle product. Take $a,b\in A$ and $u,v\in A^*$. Put
$a=f(t)\,d_qt$, $b=g(t)\,d_qt$ and
\[
   F(x) := \int_0^x (au\qsh bv) =
   \bigg(\int_0^x f(t)\int_0^t u\,d_qt\bigg)\bigg(\int_0^x g(t)\int_0^t
   v\,d_qt\bigg).
\]
Writing $F(x)=\int_0^x F'(s)\,d_qs$ where $F'(s)=(D_qF)(s)$,
and applying the $q$-product
rule for $q$-differentiation ($(D_qfg)=(D_qf)g+(\e f)(D_qg)$) yields
\begin{eqnarray*}
   F(x) &=& \int_0^x \bigg(f(s)\int_0^s u\bigg)\bigg(\int_0^s
   g(t)\int_0^t v\,d_qt\bigg)\,d_qs\\
   &&\qquad\qquad\qquad + \int_0^x g(s)\bigg(\int_0^{sq}
   f(t)\int_0^t u\,d_qt\bigg)\int_0^s v\,d_qs\\
   &=& \int_0^x \left[ a(u\qsh bv) + b(\e(au)\shuff v)\right],
\end{eqnarray*}
where the first equality follows from~(\ref{duh}) and the product
rule for $D_q$. Hence the inductive definition of the $q$-shuffle
product results. Notice that commutativity and associativity of
$\qsh$ follow immediately from~(\ref{motive}).

We will conclude with a few examples of the $q$-shuffle product
illustrating how several equivalent sums can arise from this
product. Taking $\w_1\qsh\w_2\w_3$ using the inductive definition
gives (among several possibilities):
\begin{eqnarray*}
\w_1\qsh\w_2\w_3
&=\w_1\w_2\w_3+\w_2(\e\w_1)\w_3+\w_2\w_3(\e^2\w_1)\\
&=\w_1\w_2\w_3+\w_2(\e\w_1\w_3)+\w_2\w_3(\e\w_1).
\end{eqnarray*}
Writing $\w_i=\w_i(x)$, these equations translate into the easily
verifiable generic series identities:
\begin{eqnarray*}
&&\sum_{0\le l_1}\w_1(xq^{l_1})\sum_{0\le l_2\le l_3
}\w_2(xq^{l_2})\w_3(xq^{l_3})\\ &=&\sum_{0\le l_1\le l_2\le
l_3}\w_1(xq^{l_1})\w_2(xq^{l_2})\w_3(xq^{l_3}) +\sum_{0\le l_2\le
l_1\le l_3}\w_2(xq^{l_2})\w_1(xq^{l_1+1})\w_3(xq^{l_3})\\
&&+\sum_{0\le l_2\le l_3\le
l_1}\w_2(xq^{l_2})\w_3(xq^{l_3})\w_1(xq^{l_1+2})\\ &=&\sum_{0\le
l_1\le l_2\le l_3}\!\w_1(xq^{l_1})\w_2(xq^{l_2})\w_3(xq^{l_3})
+\sum_{0\le l_2\le l_1\le
l_3}\!\w_2(xq^{l_2})\w_1(xq^{l_1+1})\w_3(xq^{l_3+1})\\
&&+\sum_{0\le l_2\le l_3\le
l_1}\w_2(xq^{l_2})\w_3(xq^{l_3})\w_1(xq^{l_1+1}).
\end{eqnarray*}
Taking the $q$-shuffle product as acting on the non-commutative
polynomials in forms, it follows that  all these different
expressions for the $q$-shuffle product tend to the ordinary
shuffle product in the limit as $q$ approaches 1.  Further results
about $q$-shuffle algebras and their combinatorics will be given
elsewhere.

\end{document}